\documentclass[12pt,leqno]{amsart}
\usepackage{amssymb,latexsym}
\usepackage{verbatim}   
\theoremstyle{plain}
\newtheorem{theorem}{Theorem}

\newtheorem{lemma}[theorem]{Lemma}
\newtheorem{proposition}[theorem]{Proposition}

\theoremstyle{definition}

\theoremstyle{remark}

\begin{document}

\newcommand{\Aut}{\textup{Aut}}
\newcommand\BigO{\textup{O}}
\newcommand{\BO}{\textup{BO}}
\newcommand{\BarHomo}{\,{}^{\overline{\hskip2.5mm}}\,}
\newcommand{\Char}{\textup{char}}
\newcommand{\CN}{\textup{C}_{\bf N}}
\newcommand{\cn}{\textup{c}_{\bf N}}
\newcommand{\CS}{\textup{C}_{\bf S}}
\newcommand{\cs}{\textup{c}_{\bf S}}
\newcommand{\diag}{\textup{diag}}
\newcommand{\E}{\mathbb{E}}
\newcommand{\End}{\textup{End}}
\newcommand{\ext}{\sqsubset}
\newcommand{\F}{\mathbb{F}}
\newcommand\GL{\textup{GL}}
\newcommand\Hom{\textup{Hom}}
\newcommand\Out{\textup{Out}}
\newcommand{\Q}{\mathbb{Q}}
\newcommand{\Symp}{\textup{Sp}}
\newcommand{\subheading}[1]{\vskip1.5mm\noindent{\sc #1}}
\newcommand\tensor{\otimes}
\newcommand\Wr{\,\textup{wr}\,}
\newcommand{\Z}{\mathbb{Z}}

\hyphenation{induced sub-fields}

\title[\tiny\upshape\rmfamily Solvable groups with
minimal composition length]{}
\date{Draft printed on \today}

\begin{center}\large\sffamily\mdseries
Solvable groups with a given solvable length,\\
and minimal composition length
\end{center}

\author{{\sffamily S.\,P. Glasby}}

\begin{abstract}
Let $\cs(d)$ denote the minimal composition
length of all finite solvable groups with
solvable (or derived) length $d$. We prove that:
\vskip3mm
\begin{center}
\begin{tabular}{|c|c|c|c|c|c|c|c|c|c|}
\hline
$d$&0&1&2&3&4&5&6&7&8\\ \hline
$\cs(d)$&0&1&2&4&5&7&8&13&15\\ \hline
\end{tabular}\;.
\end{center}
\end{abstract}

\maketitle
\centerline{\noindent 2000 Mathematics subject classification:
    20F16, 20F14, 20E34}

\section{Introduction}

\noindent
Let $\cn(d)$ (resp. $\cs(d)$) denote the minimal composition length of
a finite nilpotent group (resp. solvable group) with solvable
length~$d$. Burnside \cite{B13} knew that $\cn(0)=0$, $\cn(1)=1$,
$\cn(2)=3$, and
$\cn(3)=6$. It is shown in \cite{E00} and \cite{ENS99} that
$\cn(4)=14$. Exact values of $\cn(d)$ for $d\ge 5$ are unknown.
Hall \cite{H67,H64} showed that $2^{d-1}+d-1\le \cn(d)\le 2^d-1$. For
$d\ge 4$,
Evans-Riley improved the upper bound to $2^d-2$, and the author
(unpublished notes, 1993) improved the lower bound to $2^{d-1}+d+1$.
Mann \cite{M00} and Schneider \cite{S99b} further improved the lower
bound to $2^{d-1}+2d-4$ and $2^{d-1}+3d-10$ respectively. Upper bounds
are proved by producing specific examples. Constructing groups of
order $p^n$ and solvable length $\lfloor\log_2 n\rfloor +1$ appears
difficult, and doing so for minimal $n$ requires prescience. Such
constructions commonly do not work for small primes.

Let $\CN(d)$ (resp. $\CS(d)$) denote the set of all isomorphism
classes of finite nilpotent groups (resp. solvable groups) having
solvable length $d$, and minimal composition length. We shall
blur the distinction between a group $G$, and the isomorphism
class $[G]$ that it represents. Accordingly, we write $G\in\CN(d)$
(resp. $G\in\CS(d)$) as an abbreviation for the phrase ``$G$ is a
nilpotent group (resp. solvable group) with solvable length $d$, and
minimal composition length.'' For $G\in\CN(d)$ or $\CS(d)$,
$G^{(d-1)}$ is the unique
minimal normal subgroup of $G$ (Lemma~1(a)). [Recall that the
derived series for $G$ is defined recursively by
$G^{(0)}=G$ and $G^{(i+1)}=[G^{(i)},G^{(i)}]$ for $i\ge 0$, and the
{\it solvable} or {\it derived length} of a solvable group $G$ is the
minimal value of
$d$ such that $G^{(d)}=1$.] A major difficulty in studying groups
$G\in\CN(d)$ is that if $d>1$, then $G/G^{(d-1)}$ never lies in
$G\in\CN(d-1)$. The reason that we have made so much
progress in the solvable case is that if $G\in\CS(d)$, then
$G/G^{(k)}$ is commonly an element of $\CS(k)$ for large $k$ less than
$d$.

There is an analogy between ``minimal composition length'' groups and
$p$-groups of maximal class. The latter may be viewed as having a
given nilpotency class $c$, and minimal composition length. This class
of groups is amenable to inductive study as if $G$ has maximal class
$c$,
then $G/\gamma_c(G)$ has maximal class $c-1$. Maximal class groups
have been well studied, see \cite{H67}, III \S14.

We abbreviate the composition length of a
solvable group $G$ by $c(G)$, and its solvable (or derived) length
by $d(G)$. If $|G|=p_1^{k_1}\cdots p_s^{k_s}$, where the
$p_i$ are distinct primes, then $c(G)=k_1+\cdots+k_s$. It is clear that
\[
\cs(d)+1\le\cs(d+1)\le 2\cs(d)+1\qquad(d\ge0),
\]
where the upper bound is obtained by considering the
wreath product $G \Wr C_2$ where $G\in\CS(d)$. The above inequalities
imply that $d\le\cs(d)\le 2^d-1$. We show in the next paragraph that
$\cs(d)$ grows exponentially, and is considerably less than $\cn(d)$
for large $d$. For example, $18\le\cs(10)\le24$ and
$532\le\cn(10)\le 1022$.

If $G$ is the $r$-fold permutational wreath product $H\Wr \cdots\Wr H$
where $H=S_4$, then $|G|=|H|^{1+4+\cdots+4^{r-1}}$. Therefore
\[
c(G)=c(H)(4^r-1)/3<(4/3)\cdot4^r,\qquad\text{and}\qquad d(G)=3r.
\]
This proves that $\cs(d)<(4/3)\cdot4^{d/3}$ when $d$ is a multiple of
3. Since $9^{1/5}<4^{1/3}$, a sharper bound is obtained by taking
$H$ to be the primitive subgroup $\GL_2(3)\ltimes C_3^2$ of $S_9$. Then
$d(H)=5$ and $c(H)=7$, so
$c(G)=7(9^r-1)/8<(7/8)\cdot9^r$. Thus
$\cs(d)<(7/8)\cdot9^{d/5}$ when $d$ is a multiple of
5. Lower bounds for $\cs(d)$ require more work. It is shown in
Theorem~8 of \cite{G89} that a solvable group $G$ with $d(G)=d$ and
$c(G)=n$ satisfies
\[
d\le\alpha\log_2 n +9\qquad\text{where}\qquad \alpha=5\log_9 2+1.
\]
The smallest value of $n$ satisfying the above inequality is $\cs(d)$,
and so $2^{(d-9)/\alpha}\le\cs(d)$. Since $0.088<2^{-9/\alpha}$,
$1.3<2^{1/\alpha}$ and $9^{1/5}<1.56$, we see that
\[
(0.088)(1.3)^d<\cs(d)<(7/8)(1.56)^d\qquad\qquad(d>0)
\]
where the upper bound holds when $d$ is a multiple of 5.

In our proof that $\cs(8)\ge 15$, for example, we learn enough about
the
structure of putative groups with $d(G)=8$ and $c(G)=15$ in order to
construct them. Indeed, with more attention to detail we could
determine a complete and irredundant list of isomorphism classes
in $\CS(d)$ for $d\le8$. This requires great care
as it is all to easy to omit an isomorphism class, or to list the same
class twice. In this paper we fall short of this aim, however,
the isomorphism problem is solved for $d\le6$ in the preprint
\cite{G93}.

The groups we list in $\CS(d)$, $d\le8$, have the property that their
lattice of normal subgroups is a chain. The class of such groups,
which we call {\it normally uniserial}, is closed under quotients and
hence suited to inductive study. Moreover, if $M>N$ are normal
subgroups of a normally uniserial group, then $M/N$ is a
{\it characteristically uniserial} group, i.e. its lattice of
characteristic subgroups is a chain. Clearly, simple groups are
normally (and hence characteristically) uniserial. In \cite{G89}
the author constructs a remarkable group $G=\GL_2(3)\ltimes
3^{2+1}\ltimes
2^{6+1}\ltimes 3^{8+1}$ of order $2^{11}3^{13}$
with solvable length 10. (A more systematic construction of $G$
is given in \cite{GH92} where it is shown to be the derived 10
quotient of an infinite pro-$\{2,3\}$ group.) $G$ is normally
uniserial. I was surprised to learn that $G$ is a maximal subgroup of
the sporadic simple group Fi$_{23}$ of order
$2^{18}\cdot3^{13}\cdot5^{2}\cdot7\cdot11\cdot13\cdot17\cdot23$,
see \cite{C85}, p.~177. Indeed, $G$ has the remarkable property that
$G/G^{(d)}\in\CS(d)$ for $d=0,1,2,3,4,5,6,8$, and very likely also for
$d=10$. For the purposes of this paper it is useful to understand the
group
$G/G^{(8)}=\GL_2(3)\ltimes 3^{2+1}\ltimes 2^{6+1}$ which is described
in
\cite{N72,GH92}. Certain groups in $\CS(d)$, $d\le6$, have finite
presentations with deficiency zero, see \cite{GW93} for details.

\section{The case $d\le 6$}

\noindent
In this section we determine the solvable groups in $G\in\CS(d)$ for
$d\le6$. That is, we determine solvable groups with a given solvable
length $d\le6$, and minimal composition length subject to this
constraint. We shall determine sufficient structure of these groups
in order to compute additional values of $\cs(d)$. We stop
short of classifying the groups up to isomorphism. The determination of
$G\in\CS(d)$ for $d\le6$ is influenced by the elementary fact that a
metacyclic group is never the derived subgroup of a group. This fact
dates back to \cite{Z37}, Satz~9, p.~138.

\begin{lemma}\label{lemma1}
\begin{itemize}
\item[(a)] If $G\in\CS(d)$, then $G^{(d-1)}$ is the unique minimal
normal subgroup of $G$.
\item[(b)] Let $G$ be a solvable group with a unique minimal normal
subgroup. Let $P=\BigO_p(G)$ be nontrivial, and
suppose that $|P/\Phi(P)|$ equals $p^r$. Then $\BigO_{p'}(G)=1$, and
$G/P$
is isomorphic to a completely reducible subgroup of
$\GL_r(p)\cong\Aut(P/\Phi(P))$.
\item[(c)] If $2\le i<d(G)$, then $G^{(i-1)}/G^{(i)}$ and
$G^{(i)}/G^{(i+1)}$ are not both cyclic. In
particular, $c(G^{(i)}/G^{(i+2)})\ge3$ for $1\le i<d(G)-1$.
\item[(d)] Let $1\le i<d(G)$ and let $G^{(i-1)}/G^{(i)}$ be cyclic, and
the
unique minimal normal subgroup of $G/G^{(i)}$. Then $G/G^{(i)}$ acts
faithfully as a group of automorphisms of $G^{(i)}/G^{(i+1)}$, and
$G/G^{(i+1)}$ is a split extension of $G^{(i)}/G^{(i+1)}$ by
$G/G^{(i)}$. Moreover, $G^{(i-1)}/G^{(i)}$ has order coprime to
$|G^{(i)}/G^{(i+1)}|$ and acts fixed-point-freely.
\item[(e)] Suppose that $i\ge2$, $c(G^{(i-1)}/G^{(i)})=2$ and
$c(G^{(i)}/G^{(i+1)})=1$. Then $G^{(i-1)}/G^{(i+1)}$ is an
extraspecial group of order $p^3$.
\end{itemize}
\end{lemma}

\begin{proof}
(a) Let $N$ be a nontrivial normal subgroup of $G$. If
$G^{(d-1)}\not\le N$, then $G/N$ has solvable
length $d$, and smaller composition length. Since $G\in\CS(d)$, this is
impossible. Thus $G^{(d-1)}\le N$, as desired.
\vskip2mm\noindent
(b) The order of the unique minimal normal subgroup is a
power of some prime, say $p$, and $\BigO_{p'}(G)=1$. By a result
of Hall and Higman \cite{H67}, VI\S6.5,
$C_{G/\Phi(P)}(P/\Phi(P))=P/\Phi(P)$, and hence
\[
G/P\le\Aut(P/\Phi(P))\cong\GL_r(p).
\]
A standard argument shows that
$G/P$ acts completely reducibly, otherwise $\BigO_p(G)>P$. [Recall
that a module is called {\it completely reducible} if each submodule
has a
complementary submodule.]
\vskip2mm\noindent
(c) Suppose to the contrary that $G^{(i-1)}/G^{(i)}$ and
$G^{(i)}/G^{(i+1)}$ are both (nontrivial) cyclic groups. Then
$\Aut(G^{(i)}/G^{(i+1)})$ is abelian and so $C_G(G^{(i)}/G^{(i+1)})\le
G'$.
This implies that $G^{(i-1)}/G^{(i+1)}$ is abelian (being a cyclic
extension of a central subgroup). This is a contradiction.
\vskip2mm\noindent
(d) To simplify notation assume that $G^{(i+1)}=1$, and set
$M=G^{(i-1)}$ and $N=G^{(i)}$. Since $M/N$ is a minimal normal
subgroup of $G/N$, it is
elementary abelian. Since it is also cyclic, it has prime order,
say $p$. If $p$ divides $N$, then $M'=[M,N]<N$, a contradiction.
Thus $N$ has order coprime to $p$. Now $M\le C_G(M)<N$ because $M$ is
abelian and the chief factor $N/M$  does not
centralize $M$. Therefore, $C_G(M)=M$ and $G/M$ is a subgroup of
$\Aut(M)$. Since
$N=[M,N]\times C_N(M)$, it follows that $C_N(M)=1$, or that
$M/N$ acts fixed-point-freely on $N$. By the Frattini argument,
$G$ is a split extension of $M$ by $N_G(K)$ where $K$ is Sylow-$p$
subgroup of $M$.
\vskip2mm\noindent
(e) Since $G^{(i-1)}/G^{(i)}$ centralizes $G^{(i)}/G^{(i+1)}$, it
follows that $G^{(i-1)}/G^{(i)}$ is not cyclic. Thus there exist
primes $q$ and $p$ such that $G^{(i)}/G^{(i+1)}\cong C_q$ and
$G^{(i-1)}/G^{(i)}\cong C_p\times C_p$ . If $p\ne q$, then
$G^{(i-1)}/G^{(i+1)}$ is abelian, a contradiction. Therefore
$G^{(i-1)}/G^{(i+1)}$ is extraspecial of order $p^3$.
\end{proof}

\noindent{\bf Notation.}
Let $G$ have solvable length $d$. Write $n(G)=(n_1,n_2,\dots,n_d)$
where $n_i$ is the composition length of the abelian group
$G^{(i-1)}/G^{(i)}$. Note that $c(G)=n_1+n_2+\cdots+n_d$. The invariant
$n(G)$ will provide a useful tool for classifying elements of
$\CS(d)$. Let $K\ltimes N$ and $K\ext N$
denote a split extension, and a potentially nonsplit extension, of $N$
by $K$ respectively. Let $p$, $q$, $r$, $s$ denote primes.
Let $C_p$ and $E_p$ denote cyclic groups, and
extraspecial groups of order $p$ and $p^3$ respectively. Denote
the metacyclic group $\langle a,b\mid a^p=b^q=1,b^a=b^k\rangle$
of order $pq$ by $M_{p,q}$,
where the order of $k$ modulo $q$ is $p$. Note that $q\equiv~1\mod p$
and the isomorphism type of $M_{p,q}$ is independent of $k$.
Let $H$ denote an extension of the quaternion group of order~8 by the
symmetric group $S_3$ that has solvable length~4. There are two such
groups,
namely $\GL_2(3)$ and the binary octahedral group $\BO=\langle
a,b,c\mid
a^2=b^3=c^4=abc\rangle$.

\begin{theorem}\label{theorem2}
Let $\cs(d)$ denote minimal composition length of a finite solvable
group
with solvable length~$d$. The values of $\cs(d)$, and the structure of
$G\in\CS(d)$ for $d\le6$, are given below.
\textup{
\begin{center}
\begin{tabular}{|c|c|c|c|c|c|c|c|}
\hline
$d$&0&1&2&3&4&5&6\\ \hline
$\cs(d)$&0&1&2&4&5&7&8\\ \hline
$G$&1&$C_p$&$M_{p,q}$&$M_{p,q}\ltimes C_r^2$
&$M_{p,q}\ext E_r$&$H\ltimes C_s^2$&$H\ext E_s$
\\ \hline
&&&&$C_p\ltimes E_r$&&$\Symp_2(3)\ext E_s$&\\ \hline
\end{tabular}\;.
\end{center}
}
\end{theorem}

\begin{proof}
Let $G\in\CS(d)$. If $d\le2$, then the structure of $G$ is clear, and
hence so too are the values of $\cs(d)$. Suppose now that $d\ge3$.
It follows from Lemma~1(c) that $n_i+n_{i+1}\ge 3$ for $i\ge2$.
Hence the possible values of $n(G)$ are:
\begin{center}
\begin{tabular}{|c|c|c|c|c|c|c|c|}
\hline
$d$&3&4&5&6\\ \hline
$n(G)$&$(1,1,2)$&$(1,1,2,1)$&$(1,1,2,1,2)$&$(1,1,2,1,2,1)$
\\ \hline
&$(1,2,1)$&&$(1,2,1,2,1)$&\\ \hline
\end{tabular}\;.
\end{center}

The question arises as to whether each of the 6 above values of $n(G)$
arise for particular groups $G$. The answer is affirmative.
There is a subgroup of the automorphism group of an exponent-$p$
extraspecial group of order $p^{2k+1}$ isomorphic to the general
symplectic group $\textup{GSp}_{2k}(p)$, see \cite{W72,GH92}.
Thus we may form the split extension $\textup{GSp}_{2k}(p)\ltimes
p^{2k+1}$. When $k=1$ and $p=3$ this group is $G=\GL_2(3)\ltimes E_3$
as
$\textup{GSp}_2(3)\cong\GL_2(3)$. Now $G$ has solvable length~6, and
the
quotients $G^{(i-1)}/G^{(i)}$ are $C_2$, $C_3$, $C_2\times C_2$,
$C_2$, $C_3\times C_3$, $C_3$. Thus $n(G)$ equals $(1,1,2,1,2,1)$. By
taking
quotients of $G$ or $G'$ we see that each of the 6 above values
of $n(G)$ arise.

We shall now be more specific about the structure of an arbitrary
group $G$ such that $n(G)$ is one of the 6 above values. It is clear
that $G\in\CS(d)$. If $n(G)=(1,1,2)$, then $G^{(2)}$ is not cyclic by
Lemma~1(c), so $G^{(2)}\cong C_r^2$ for some prime $r$. By Lemma~1(d),
$G$ is a split extension $M_{p,q}\ltimes C_r^2$. Indeed,
$M_{p,q}\le\GL_2(r)$ acts irreducibly. If $n(G)=(1,2,1)$, then
$G'=E_r$ is extraspecial of order $r^3$ by Lemma~1(e), and
$G=C_p\ltimes E_r$ where $C_p$ acts fixed-point-freely on
$E_r/\Phi(E_r)$. When $n(G)=(1,1,2,1)$, then $G=M_{p,q}\ext E_r$.
Suppose that $n(G)=(1,1,2,1,2)$. Then $G^{(4)}$ is noncyclic, say
$C_s^2$ where $s$ is prime. Now $G^{(2)}/G^{(4)}$ is an extraspecial
group by Lemma~1(e), and it acts irreducibly on $G^{(4)}$. This forces
$G^{(2)}/G^{(4)}$ to be isomorphic to the quaternion group $Q_8$, or
the dihedral group $D_8$, of order~8. As $\Out(D_8)\cong C_2$, and
$\Out(Q_8)\cong S_3$, it follows that $H=G/G^{(4)}$ is an extension of
$Q_8$ by $S_3$. Therefore, $H\cong\GL_2(3)$ or $\BO$. By Lemma~1(d),
$G$ is a split extension $H\ltimes C_s^2$. The action of $H$ on $C_s^2$
 is
irreducible, and exists only for certain odd primes $s$.
Arguing as above, the structure of $G$ satisfying
$n(G)=(1,2,1,2,1)$ is $\Symp_2(3)\ext E_s$, where $\Symp_2(3)$ denotes
the symplectic group and $\Symp_2(3)\cong H'$. If $s=3$, then there are
nonisomorphic split and nonsplit extensions of $E_s$ by $\Symp_2(3)$,
see \cite{G93}.
Finally when $n(G)=(1,1,2,1,2,1)$, $H=G/G^{(4)}\cong\GL_2(3)$ or $\BO$
and $G^{(4)}\cong E_s$ is extraspecial of order $s^3$.
\end{proof}

\section{The case $d=7$}

\noindent
Before proving that $\cs(7)=15$ in Theorem~7, we need 4 preliminary
lemmas.

\begin{lemma}\label{lemma3}
Let $\textup{cr}(n)$ denote the maximal solvable
length of a completely reducible solvable subgroup of $\GL_n(\F)$,
where the field $\F$ may vary. Then
\textup{
\begin{center}
\begin{tabular}{|c|c|c|c|c|c|c|c|c|}
\hline
$n$&1&2&3&4&5&6&7&8\\ \hline
$\textup{cr}(n)$&1&4&5&5&5&6&6&8\\ \hline
\end{tabular}\;.
\end{center}
}
\end{lemma}

\begin{proof} See \cite{N72} for an explicit formula for the function
$\textup{cr}(n)$.
\end{proof}

\begin{lemma}\label{lemma4}
Let $P$ be a finite abelian group, and let $Q$ be a
solvable subgroup of $\Aut(P)$ with solvable length~$d$.
\begin{itemize}
\item[(a)] Then
\textup{
\begin{center}
\begin{tabular}{|c|c|c|c|c|}
\hline
$|P|$&$p$&$p^2$&$p^3$&$p^4$\\ \hline
$d$&$\le1$&$\le4$&$\le5$&$\le6$\\ \hline
\end{tabular}\;.
\end{center}
}
\item[(b)] A subgroup chain $P=P_0>P_1>\cdots>P_n=1$ is called
\emph{maximal} if $|P|=p^n$, or equivalently $|P_{i-1}:P_i|=p$ for
$i=1,\dots,n$. If $P$ is an abelian group of order dividing $p^4$,
and $Q$ stabilizes a maximal subgroup chain, then $d\le3$.
\end{itemize}
\end{lemma}

\begin{proof} (a) If $P$ is elementary abelian of order $p^n$, then
the maximum value of $d$ is given in \cite{N72}, Theorem~A. In
particular,
$d=1,4,5,6$ when $n=1,2,3,4$. If $P$ is not elementary abelian, then
$P^p=\{g^p\mid g\in P\}$ is a proper nontrivial characteristic
subgroup. Furthermore, the automorphisms of $P$ centralizing both
$P/P^p$ and $P^p$, form an abelian group. The above table follows from
these two facts.
\vskip2mm\noindent
(b) This follows when $P$ is elementary abelian, as then $Q$
is a subgroup of the upper triangular matrices. If $P$ is not
elementary abelian, then consider the groups $P/P^p$ and $P^p$ as
above.
\end{proof}

Much more is known about primitive maximal solvable linear groups than
is given in the following lemma, however, this simplified form is all
that we require.

\begin{lemma}\label{lemma5}
Let $M$ be an absolutely irreducible primitive maximal solvable
subgroup of $\GL_r(\F)$ where $\F$ is a finite field. Then $Z:=Z(M)$
is cyclic of order $|\F|-1$. If $F$ is the Fitting radical of $M$ (the
maximal nilpotent normal subgroup of $M$), then $F/Z$ is elementary
abelian of order $r^2$. If $r=p_1^{k_1}\cdots p_s^{k_s}$ where the
$p_i$ are distinct primes, then there exist extraspecial subgroups
$E_i$ of $F$ of order $p_i^{2k_i+1}$ such that $F$ is a central
product $(E_1\times\cdots\times E_s)\textup{\sffamily Y} Z$, and $F$ is
conjugate in $\GL_r(\F)$ to $(E_1\tensor\cdots\tensor E_s)Z$.
\end{lemma}

\begin{proof}
The first two sentences follow from \cite{S76}, Lemma~19.1 and
Theorem~20.9, and the last sentence can be deduced from results on
pages
141--146. A more convenient reference is \cite{S92}, Theorems~2.5.13
and 2.5.19.
\end{proof}

The following result is proved in \cite{ENS99,S99a}.

\begin{lemma}\label{lemma6}
Let $p\ge3$ be a prime, and let $P$ be a $p$-group satisfying
$|P'/P''|=p^3$ and $P''\ne1$. Then
\[
P'=\gamma_2(P)>\gamma_3(P)>\gamma_4(P)>\gamma_5(P)=P''.
\]
\end{lemma}

\begin{theorem}\label{theorem7}
A finite solvable group with solvable length $7$ has
composition length at least $13$, and this bound is best possible. More
succinctly, $\cs(7)=13$.
\end{theorem}

\begin{proof}
As usual, our proof has two parts: (1) show that if $d(G)=7$, then
$c(G)\ge13$, and (2) exhibit a group $G$ with $d(G)=7$ and $c(G)=13$.
The second part is deferred to Proposition~8 below.

Suppose that $d(G)=7$. By the proof of Lemma~1(a), we may reduce to
the case that $G$ has a unique minimal normal subgroup. By Lemma~1(b),
there is a (unique) prime $p$ such that $P:=\BigO_p(G)$ is nontrivial,
and
$Q:=G/P$ is a completely reducible subgroup of $\GL_r(p)$ where
$|P/\Phi(P)|=p^r$. If $d(P)\ge 4$, then $c(P)\ge 14$
(\cite{E00,ENS99}),
and hence $c(G)\ge 14$. If $d(P)=1$, then $d(Q)\ge 6$ and
$c(Q)\ge\cs(6)=8$ by Theorem~2. However, $r\ge6$ by Lemma~3, and so
\[
 c(G)=c(Q)+c(P)\ge 8+6=14.
\]
The two remaining cases when $d(P)=2$ or 3
require more detailed analyses.

\subheading{Case (a) $d(P)=2$.}
Now $d(Q)\ge5$, so $c(Q)\ge\cs(5)=7$ by Theorem~2. If $c(P)\ge6$, then
$c(G)=c(Q)+c(P)\ge 7+6=13$. Thus it suffices to consider the cases
when $c(P)<6$. Let $|P|=p^{r+s}$ where $|\Phi(P)|=p^s$. Then $r\ge3$ by
Lemma~3 and there are three cases when $d(P)<6$, namely
\[
(r,s)=(3,1), (3,2)\quad\text{and}\quad (4,1).
\]
We show that the first possibility never arises, and if the second or
third arise, then $c(G)\ge13$.

\subheading{Subcase $(r,s)=(3,1)$.} In this case $\Phi(P)=P'$ has
order~$p$. If $Z(P)=P'$, then $P$ is an extraspecial group with even
composition length, a contradiction. Hence $\Phi(P)<Z(P)<P$ and since
$Q$
acts completely reducibly, $Q\le\GL_1(p)\times\GL_2(p)$ by
Lemma~1(b). Thus $d(Q)\le4$ by
Lemma~3. This is a contradiction as $d(Q)\ge5$. Hence this case never
arises.

\subheading{Subcase $(r,s)=(3,2)$.} Arguing as in the previous case,
we see that $Q\le\GL_3(p)$ is an irreducible subgroup. If $Q$ does not
act absolutely irreducibly, then $Q\le\GL_1(p^3)$, and $d(Q)\le1$, a
contradiction. If $Q\le\GL_3(p)$ is an imprimitive subgroup, then
$Q\le\GL_1(p)\Wr S_3$ and $d(Q)\le3$, a contradiction. In summary,
$Q\le\GL_3(p)$ acts absolutely irreducibly and primitively. Thus $Q$
is a subgroup of an absolutely irreducible primitive maximal solvable
subgroup $M$ of
$\GL_3(p)$. By Lemma~5 there are characteristic subgroups $Z\le F\le
M$ such that $F/Z$ is elementary abelian of order $3^2$,
$M/F\le\Symp_2(3)$, and $F'$ has order 3. Since $c(Q)\ge7$ and
$c(P)=5$,
we must eliminate the case when $c(G)=12$. In this case, $c(Q)=7$,
$M/F\cong\Symp_2(3)$, and $F$ contains an extraspecial subgroup
of order $3^3$ and exponent 3, and $M=QZ$. Since $Z(Q)\le Z(M)$,
and $M$ acts absolutely irreducibly, there
is an element $z\in Z(Q)$ of order~3 which induces the scalar
transformation
$\omega1$ on $P/\Phi(P)$ where $\omega$ is primitive cube root
of~1 modulo~$p$. We view $z$ as an element of $G^{(4)}$ of order 3.

We show that $\Phi(P)\le Z(P)$. If $\Phi(P)\not\le Z(P)$, then
\[
\Phi(P)<Z(P)\Phi(P)<P.
\]
This contradicts the fact that $Q$ acts irreducibly on $P/\Phi(P)$.
In summary, we know that $\Phi(P)\le Z(P)$, and
$g^z\Phi(P)=g^\omega\Phi(P)$ for all $g\in P$. Therefore,
\[
[g_1,g_2]^z=[g_1^\omega,g_2^\omega]=[g_1,g_2]^{\omega^2}
\qquad(g_2,g_2\in P).
\]
This proves that $z$ acts nontrivially on $P'$, and hence $\Aut(P')$
contains a subgroup with solvable length at least 5, contrary to
Lemma~4(a). Thus we have proved $c(G)\ge13$ in this case.

\subheading{Subcase $(r,s)=(4,1)$.}
Then $\Phi(P)=P'$ has order $p$, so $P'\le Z(P)$.
Since $d(P)>1$, it follows that $p^2\le|P:Z(P)|\le p^4$. If
$|P:Z(P)|=p^2$, then it follows from Lemma~1(b) that
$Q\le\GL_2(p)\times\GL_2(p)$, and hence $d(Q)\le4$ by Lemma~3. This
is a contradiction as $d(Q)\ge5$. If $|P:Z(P)|=p^3$, then similar
reasoning shows $Q\le\GL_3(p)\times\GL_1(p)$. Arguing as in the
previous subcase, $Q$ acts absolutely irreducibly and primitively on
$P/Z(P)$. Appealing as above to Lemma~5, a nontrivial element
$z\in Z(Q)$ maps generators $a_1,a_2,a_3,a_4$ for $P$ to
$a_1^\omega,a_2^\omega,a_3^\omega,a_4$ modulo $\Phi(P)$ where
$\omega$ is a primitive cube root of unity modulo~$p$.
Since $P'=\Phi(P)\le Z(P)$, there is a well defined action of $z$ on
$P'$. Since $[a_i,a_j]^z=[a_i^z,a_j^z]=[a_i,a_j]^\omega$ or
$[a_i,a_j]^{\omega^2}$, $z\in G^{(4)}$ acts nontrivially on $P'$. Thus
$\Aut(P')$ contains a subgroup with
solvable length at least 5, contrary to Lemma~4. Thus $c(G)\ge13$ in
this case also.

\subheading{Case (b) $d(P)=3$.}
Since $P''\ne1$, $|P'/P''|\ge p^3$ by Hilfsatz~7.10 of \cite{H67}.
If $p=2$, then $\GL_2(2)$ and $\GL_3(2)$ are too small to accommodate
a solvable subgroup $Q$ with $d(Q)\ge4$ and $c(Q)\ge\cs(4)=5$. Hence
if $p=2$, then $r\ge4$ and
\[
c(P)\ge c(P/\Phi(P))+c(P'/P'')+c(P'')\ge 4+3+1=8.
\]
Therefore $c(G)=c(Q)+c(P)\ge 5+8=13$ as desired. Assume now that
$p\ge3$ and $|P'/P''|= p^3$. By Lemma~6,
\[
P'=\gamma_2(P)>\gamma_3(P)>\gamma_4(P)>\gamma_5(P)=P''.
\]
Now $P/P'$ acts nontrivially on $P'/P''$. Since $G^{(3)}\not\le P$,
it follows that
$\Aut(P'/P'')$ contains a subgroup with solvable length at least~4.
This contradicts Lemma~4(b). Henceforth assume that
$|P'/P''|\ge p^4$.

In summary, $c(Q)\ge5$ and $c(P)\ge7$, so $c(G)\ge12$. Assume by way
of contradiction that $c(G)=12$. Then $c(Q)=5$ and $c(P)=7$. Since
$d(Q)=4$, we have $Q\in\CS(4)$. Thus $Q\cong\GL_2(3)$ or $\BO$ by
Theorem~2. In addition, $|P/P'|=p^2, |P'/P''|=p^4$ and $|P''|=p$.
Thus $\gamma_2(P)/\gamma_3(P)$ is cyclic, and so
\[
P''=[\gamma_2(P),\gamma_2(P)]=[\gamma_2(P),\gamma_3(P)]\le\gamma_5(P).
\]
If $P''<\gamma_5(P)$, then
$P'=\gamma_2(P)>\gamma_3(P)>\gamma_4(P)>\gamma_5(P)>P''$. However,
$P/P'$ acts nontrivially on the abelian group $P'/P''$ of order~$p^4$.
As $Q$ acts irreducibly on $P/P'$, it follows that $G^{(4)}=P$.
Thus $\Aut(P'/P'')$ contains a subgroup with solvable
length at least~4, contrary to Lemma~4(b). Hence $P''=\gamma_5(P)$.

If the cyclic group $\gamma_2(P)/\gamma_3(P)$ has order at least
$p^2$, then its order is exactly $p^2$, and we have the characteristic
series
\[
P'=\gamma_2(P)>\gamma_2(P)^p\gamma_3(P)>\gamma_3(P)>
\gamma_4(P)>\gamma_5(P)=P''.
\]
As above, this is impossible. Thus $|\gamma_2(P)/\gamma_3(P)|=p$, and
$|\gamma_3(P)|=p^4$. Now $\gamma_3(P)$ is abelian as
$[\gamma_3(P),\gamma_3(P)]\le\gamma_6(P)=1$. Exactly one of
$|\gamma_3(P):\gamma_4(P)|$ or $|\gamma_4(P):\gamma_5(P)|$ has order
$p^2$. Suppose that $\gamma_3(P)$ has a characteristic subgroup $N$ of
index $p^2$, and $K$ is a solvable group of automorphisms of
$\gamma_3(P)$. By Lemma~4(a), $K^{(4)}$ centralizes both
$\gamma_3(P)/N$ and
$N$. Since $N$ is abelian, it follows that $K^{(5)}=1$.
However, $P'/\gamma_3(P)$ acts nontrivially on $\gamma_3(P)$ and
$G^{(4)}\not\le P'$, so $\Aut(\gamma_3(P))$ contains a subgroup with
solvable length at least~6. This contradicts the fact that $K^{(5)}=1$,
and proves
that $|\gamma_4(P):\gamma_5(P)|=p^2$. Now $K=G/\gamma_3(P)$ satisfies
$d(K)=6$ and $c(K)=8$. Thus $K\in\CS(6)$. By Theorem~2, $K\cong
H\ext E_p$ where $H\cong\GL_2(3)$ or $\BO$.

Consider the section $G^{(3)}/\gamma_4(P)$. Since $G^{(4)}=P$ we have
\[
|G^{(3)}:P|=2, |P:P'|=p^2\quad\text{and}\quad
|P':\gamma_3(P)|=|\gamma_3(P):\gamma_4(P)|=p.
\]
Let $z\in G^{(3)}$ have order 2. It follows from the
structure of $G/\gamma_3(P)$ that $z$ acts as the scalar transformation
$-I$
on $P/P'\cong C_p\times C_p$. As $p$ is odd, and $z$ centralizes both
$P'/\gamma_3(P)$ and $\gamma_3(P)/\gamma_4(P)$, it centralizes
the abelian group $P'/\gamma_4(P)$. Let $g\in P$ and $h\in
P'$. Then
\[
[g,h]\equiv [g,h]^z\equiv [g^z,h^z]\equiv [g^{-1},h]
\equiv [g,h]^{-1}\,\text{mod}\,\gamma_4(P).
\]
As $p$ is odd and $[g,h]^2\equiv 1\mod\gamma_4(P)$, we see
that $[P,P']\subseteq\gamma_4(P)$. This is a contradiction as
$\gamma_3(P)\not\subseteq\gamma_4(P)$.

In summary, we have proved in each case that if $d(G)=7$, then
$c(G)\ge 13$.
\end{proof}

The last case in Theorem~7 was difficult to eliminate.
We can show that $\gamma_3(P)$ is either $(C_p)^4$ or
$C_{p^2}\times(C_p)^2$. In either case, there exist a subgroup
$H\ltimes E_p$ of $\Aut(\gamma_3(P))$ with solvable length~6
normalizing
a subgroup chain $\gamma_3(P)=P_0>P_1>P_2>P_3=1$ with
$|P_0:P_1|=|P_2:P_3|=p$, and $|P_1:P_2|=p^2$. Our contradiction was
therefore subtle. It arose not because the action of
$G/\gamma_3(P)$ on $\gamma_3(P)$ was untenable, rather because
there was no extension of
$\gamma_3(P)$ by $H\ltimes E_p$ having solvable length 7.

\begin{proposition}\label{proposition8}
There exists a solvable group with solvable length~$7$ and
composition length~$13$. Thus $\cs(7)\le13$.
\end{proposition}

\begin{proof}
Let $V$ be an $r$-dimensional vector space over a field $\F$.
The homogeneous component $\Lambda^i V$ of the exterior algebra
$\oplus_{i=0}^r \Lambda^i V$ has dimension ${r\choose i}$. Set
$P=V\times\Lambda^2 V$,
and define a binary operation on $P$ via the rule
\[
(v_1,w_1)(v_2,w_2)=(v_1+v_2,w_1+w_2+v_1\wedge v_2)
\]
where $v_1,v_2\in V,w_1,w_2\in \Lambda^2V$.
Then $P$ is a group. If $\Char(\F)\ne2$, then the derived subgroup
$P'$ equals $0\times\Lambda^2V$ because $v_2\wedge v_1\ne-v_1\wedge
v_2$.
The right action of $\GL_r(\F)$ on $P$ defined by
$(v,w)g=(vg,w(g\wedge g))$ gives rise to a split extension
$GL_r(\F)\ltimes P$. We are interested in the subgroup
$K\ltimes P$ of this group when $r=3$, $|\F|=p$ is an odd prime and
$K\le\GL_3(p)$ is isomorphic to $\Symp_2(3)\ltimes E_3$.
If $p\equiv1\mod~3$, then there are faithful representations
$\Symp_2(3)\ltimes E_3 \to\GL_3(p)$. [Indeed, when $p\equiv1\mod~9$,
then there are faithful representations of the nonsplit extensions
$\Symp_2(3)\ext E_3 \to\GL_3(p)$.] Let $G=K\ltimes P$. Then $c(K)=7$
and $c(P)={3\choose 1}+{3\choose 2}=6$, so $c(G)=c(K)+c(P)=13$. We
show now that $d(G)=7$. An element $z\in K^{(4)}$ of order 3
induces the scalar transformation $\omega 1$ on $P/\Phi(P)\cong V$,
where
$\omega$ has order 3 modulo~$p$.
If $k\in K$ has matrix $A$ relative to a basis $e_1,e_2,e_3$ for $V$,
then $k\wedge k$ has matrix $\text{det}(A)(A^{-1})^T$ relative to the
basis
$e_2\wedge e_3,e_3\wedge e_1,e_1\wedge e_2$ for
$\Lambda^2V$. Therefore, $z$ acts like $\omega^2 1$ on $\Phi(P)=P'$.
This shows that $G^{(5)}=P$, and hence that $d(G)=7$.
\end{proof}

\section{The case $d=8$}

\noindent

\begin{theorem}\label{theorem9}
A finite solvable group with solvable length $8$ has
composition length at least $15$, and this bound is best possible. More
succinctly, $\cs(8)=15$.
\end{theorem}

\begin{proof}
As remarked in the introduction, the  group $\GL_2(3)\ltimes
E_3\ltimes 2^{6+1}$ of order $2^{11}3^4$ has solvable length~8.
This proves that $\cs(8)\le 15$. Since $\cs(7)=13$, we see that
$\cs(8)=14$ or 15. We eliminate the case $\cs(8)=14$.

Let $G\in\CS(8)$. Suppose that $P=\BigO_p(G)$ is nontrivial and
$|P/\Phi(P)|$ equals $p^r$. Then
$Q=G/P$ is a completely reducible subgroup of $\GL_r(p)$.
If $d(P)\ge 4$, then $c(P)\ge14$ by \cite{ENS99},
and hence $c(G)>15$. If $d(P)=1$, then $d(Q)\ge7$ and $r\ge8$ by
Lemma~3. By Theorem~7, $c(Q)\ge13$ so
$c(G)\ge 13+8=21$. We shall now consider the two remaining cases:
$d(P)=2$ or 3.

\subheading{Case $d(P)=2$.}
Now $d(Q)\ge6$, so $c(Q)\ge\cs(6)=8$. By Lemma~3, $r\ge6$
therefore $c(P)\ge7$, and so $c(G)\ge8+7=15$.

\subheading{Case $d(P)=3$.}
Now $d(Q)\ge5$, so $c(Q)\ge\cs(5)=7$. By Lemma~3, $r\ge3$. Since
$|P'/P''|\ge p^3$, it follows that $|P|\ge p^7$. Therefore
$c(G)\ge7+7=14$.
Suppose that $c(G)=14$. Then $c(Q)=7$, $|P|=p^7$, $P'=\Phi(P)$,
$|P':P''|=p^3$ and $|P''|=p$. By Lemma~3, $Q\le\GL_3(p)$ acts
irreducibly. Arguing as in Theorem~2, $Q$ acts absolutely irreducibly
and
primitively. Therefore, $G^{(5)}=P$. It follows
from Lemma~5 that $Q\cong\Symp_2(3)\ext E_3$ where $E_3$ has exponent
3. Now $P/P'$ acts nontrivially on $P'/P''$. Therefore $\Aut(P'/P'')$
contains a subgroup with solvable length at least 6. This is
impossible by Lemma~4(b). Hence $\cs(8)=15$ as claimed.
\end{proof}

With more precise arguments, we can show that if $d(G)=8$ and
$c(G)=15$, then $d(P)=2$. By Lemma~3, $G/P$ acts irreducibly on
$P/\Phi(P)$, and so $\Phi(P)=P'=Z(P)$. Thus $P$ is an extraspecial
group of order $p^{6+1}$ (and exponent $p$, if $p$ is odd). Since
$c(Q)=8$ and $d(Q)=6$, $Q\cong H\ext E_s$ by Theorem~2. The
representation theory of extraspecial groups shows that $s=3$, and
hence $Q\cong\GL_2(3)\ltimes E_3$. In addition, $p\equiv -1\mod~3$, and
$Q'$ acts irreducibly but not absolutely irreducibly on
$P/\Phi(P)$. In summary, elements of $\CS(8)$ have the form
$(\GL_2(3)\ltimes E_3)\ext p^{6+1}$.

\vskip5mm
\begin{center}
{\sc Acknowledgements}
\end{center}
\vskip2mm

\noindent
{\small
I would like to thank C.W. Parker for alerting me to the Fi$_{23}$
connection, and R.B. Howlett for discussions which led to
the examples in Proposition~8.
}

\vskip3mm
\goodbreak
\def\efont{\scriptsize\upshape\ttfamily}
{\tiny\scshape
\begin{tabbing}
\=\hspace{70mm}\=\kill\\
\>Department of Mathematics    \\
\>Central Washington University\\
\>WA 98926-7424, USA           \\
\>\efont GlasbyS@cwu.edu       \\
\end{tabbing}
}

\end{document}